\documentclass[twoside]{article}
\usepackage{amssymb}
\usepackage{indentfirst}
\usepackage{dsfont}
\usepackage{ntheorem}
\usepackage{amsmath}

\begin{document}
\baselineskip 13pt \noindent \thispagestyle{empty}

\markboth{\centerline{\rm Z. Tu \; \& \; L. Wang}}{\centerline{\rm
Classification of proper self-mappings}}

\begin{center} \Large{\bf  Classification of proper holomorphic mappings between
certain unbounded non-hyperbolic domains }
\end{center}

\begin{center}
\noindent\text{Zhenhan Tu$^1$,  \; Lei Wang$^2$$^*$ }\\
\vskip 4pt \noindent\small {$^1$School of Mathematics and Statistics,
Wuhan University, Wuhan, Hubei 430072, P.R. China}\\
\vskip 4pt \noindent\small {$^2$ School of Mathematics and Statistics, Huazhong University of Science and Technology, \\ Wuhan, Hubei 430074, P.R. China}\\
\vskip 4pt \noindent\text{Email: zhhtu.math@whu.edu.cn (Z. Tu),\; wanglei2017@hust.edu.cn (L. Wang)}
\renewcommand{\thefootnote}{{}}
\footnote{\hskip -16pt {$^{*}$Corresponding author. \\ } }
\end{center}

\begin{center}
\begin{minipage}{13cm}
{\bf Abstract.} {\small The Fock-Bargmann-Hartogs domain $D_{n,m}(\mu)$ ($\mu>0$)
 in $\mathbb{C}^{n+m}$ is defined by the inequality $\|w\|^2<e^{-\mu\|z\|^2},$
 where $(z,w)\in \mathbb{C}^n\times \mathbb{C}^m$, which is an unbounded non-hyperbolic domain
 in $\mathbb{C}^{n+m}$. Recently, Tu-Wang obtained the rigidity result
that proper holomorphic self-mappings of $D_{n,m}(\mu)$ are automorphisms for $m\geq 2$,
 and found a counter-example to show that the rigidity result isn't true for $D_{n,1}(\mu)$.
 In this article, we obtain a classification of proper holomorphic mappings between $D_{n,1}(\mu)$ and $D_{N,1}(\mu)$ with $N<2n$.

 \vskip 5pt
 {\bf Key words:} Fock-Bargmann-Hartogs domains, local biholomorphisms, proper holomorphic mappings
\vskip 5pt
 {\bf 2010 Mathematics Subject Classification:} Primary 32A07,\; 32H35,\; 32M05. }
\end{minipage}
\end{center}

\numberwithin{equation}{section}
\def\theequation{\arabic{section}.\arabic{equation}}

\section{Introduction}
In 1907, Poincar\'{e} \cite{Poin} proved the following result for $n=2$:
\newtheorem{lem}{Lemma}[section]
\newtheorem{thm}{Theorem}[section]
\newtheorem{prop}{Proposition}[section]
\newtheorem{cor}{Corollary}[section]

\vskip 6pt \noindent {\bf Theorem 1.A}  {\it
Let $\mathbb{B}^n=\{z\in\mathbb{C}^n: |z|<1\}$ be the unit ball and let $U=\{|z-a|<\varepsilon\}$ be a neighborhood of a boundary point $a$ of $\mathbb{B}^n$. If $f:\mathbb{B}^n\cap U\rightarrow \mathbb{C}^n$ is a biholomorphic mapping such that  $ f\in C^1(\overline{\mathbb{B}^n}\cap U)$ and $f(\partial \mathbb{B}^n\cap U)\subset \partial \mathbb{B}^n$, then for $n>1$ this mapping extends to a biholomorphic automorphism of the whole ball,
and hence, is linear fractional.}

\vskip 6pt
This was proved for arbitrary $n>1$ in 1962 by Tanaka \cite{Tana} who was
apparently not aware of Poincar\'{e}'s work. The same result was rediscovered
by Alexander \cite{Ale} and Pelles \cite{Pel}. Pin\v{c}uk \cite{Pin75, Pin78} established a new proof of the previous result and extended it to strongly
pseudoconvex domains in $\mathbf{C}^n$  with real analytic simply connected boundaries.

When $n=1$, Theorem 1.A is no long true. However, the following classification is obvious.

\vskip 6pt
\noindent
{\bf Theorem 1.B}  {\it Every proper holomorphic mapping $f$ from the unit disc to itself is a finite Blaschke product. That is, there are finitely many points $\{a_j\}$ in the disc, positive integer multiplicities $m_j$, and a number $e^{i\theta}$ such that
$$f(z)=e^{i\theta} \prod^m_{j=1} \Big(\frac{a_j-z}{1-\overline{a_j} z} \Big)^{m_j}.$$ }

Alexander \cite{Ale, Ale2} further studied
proper holomorphic mappings between bounded domains with the same dimension. Alexander's theorem has been generalized to many classes of domains (e.g., see Bedford-Bell \cite{BB},  Diederich-Forn{\ae}ss \cite{DF}, Huang  \cite{H},  Su-Tu-Wang \cite{ST},  Tu  \cite{Tu},  Tu-Wang  \cite{TW,TW2}, and Webster \cite{Web}). Inspired by these theorems, there are many results on classifying proper holomorphic mappings up to holomorphic automorphisms (e.g., see Dini-Primicerio  \cite{DS}, Ebenfelt-Son  \cite{Eben}, Faran  \cite{Far}, Landucci-Pinchuk  \cite{LP}, Spiro \cite{Spi}, and Zapalowski \cite{Zapa}).

\vskip 6pt
The Fock-Bargmann-Hartogs domains $D_{n,m}(\mu)$ are defined by
$$D_{n,m}(\mu):=\{(z,w)\in\mathbb{C}^n\times\mathbb{C}^m:\|w\|^2<e^{-\mu\|z\|^2}\}, \;\;\; \mu>0.$$
The Fock-Bargmann-Hartogs domains $D_{n,m}(\mu)$ are unbounded
strongly pseudoconvex domains in $\mathbb{C}^{n+m}$. We note that
each $D_{n,m}(\mu)$ contains $\{(z, 0)\in
\mathbb{C}^n\times\mathbb{C}^m\} \cong \mathbb{C}^n$. Thus each
$D_{n,m}(\mu)$ is not hyperbolic in the sense of Kobayashi and
$D_{n,m}(\mu)$ can not be biholomorphic to any bounded domain in
$\mathbb{C}^{n+m}$. Therefore, each Fock-Bargmann-Hartogs domain
$D_{n,m}(\mu)$ is an unbounded non-hyperbolic domain in
$\mathbb{C}^{n+m}.$

In 2014, by checking that the Bergman
kernel ensures revised Cartan's theorem, Kim-Thu-Yamamori \cite{Kim}
determined the holomorphic automorphism groups of the Fock-Bargmann-Hartogs
domains as follows:

\vskip 6pt

\noindent {\bf Theorem 1.C} (Kim-Thu-Yamamori \cite{Kim}) {\it The
automorphism group ${\rm Aut}(D_{n,m}(\mu))$ is exactly the group
generated by the following automorphisms of  $D_{n,m}(\mu)$:
\begin{equation*}
\begin{array}{l}
\varphi_{U}:(z,w)\longmapsto(Uz,w),\quad U\in \mathcal{U}(n); \\
\varphi_{U^{'}}:(z,w)\longmapsto(z,U^{'}w),\quad U^{'}\in \mathcal{U}(m);\\
\varphi_{v}:(z,w)\longmapsto(z+v,e^{-\mu \langle z,v
\rangle-\frac{\mu}{2} {\left\lVert v\right\rVert}^{2}}w),\quad (v\in
\mathbb{C}^{n}),
\end{array}
\end{equation*}
where $\mathcal{U}(k)$ is the unitary group of degree $k,$ and
$\langle \cdot,\cdot \rangle$ is the standard Hermitian inner product on
$\mathbb{C}^{n}$.}

\vskip 6pt
\noindent
{\bf Remark}. According to Theorem 1.C, although the group ${\rm Aut}(D_{n,m}(\mu))$ is not transitive on the domain itself as there is no way to map a point $(z,w)$ where $w\neq 0$ to $(z',0)$ as none of the generators allow that, ${\rm Aut}(D_{n,m}(\mu))$ is transitive on the boundary $\partial D_{n,m}(\mu)$. This is a key point for the proof of our results in this paper.

\vskip 6pt
Recently, Tu-Wang \cite{TW} obtained the rigidity result on proper holomorphic mappings between two equidimensional Fock-Bargmann-Hartogs domains.

\vskip 6pt \noindent {\bf Theorem 1.D} (Tu-Wang \cite{TW})  {\it
If $D_{n,m}(\mu)$ and $D_{n',m'}(\mu')$ are two equidimensional Fock-Bargmann-Hartogs domains with $m\geq 2$
and $f$ is a proper holomorphic mapping from $D_{n,m}(\mu)$ into $D_{n',m'}(\mu')$, then $f$ is a biholomorphism between $D_{n,m}(\mu)$ and $D_{n',m'}(\mu')$.}

\vskip 6pt
\noindent {\bf For example.}  Let
$$\Phi(z_1,\cdots,z_n,w_1):=(\sqrt{2}z_1,\cdots,\sqrt{2}z_n,w_1^2),
\quad (z_1,\cdots,z_n,w_1)\in D_{n,1}(\mu).$$
Then $\Phi$ is a proper
holomorphic self-mapping of $D_{n,1}(\mu)$, but it is  branched and
isn't an automorphism of $D_{n,1}(\mu)$. Thus the assumption ``$m\geq
2$" in Theorem 1.1 cannot be removed. Also, this example implies
that a proper holomorphic self-mapping of unbounded strongly
pseudoconvex domain in $\mathbb{C}^n(n\geq 2)$ is possibly not an
automorphism. This seems to be much different from the bounded case (cf. Bedford-Bell \cite{BB} and Diederich-Forn{\ae}ss \cite{DF}).

In this article, we obtain a classification of proper holomorphic self-mappings of $D_{n,1}(\mu)$ as follows:
\begin{thm}
  Let $D_{n,1}(\mu)$ be a Fock-Bargmann-Hartogs domain and $F$ be a proper holomorphic self-mapping of $D_{n,1}(\mu)$.
  Then there exist holomorphic automorphisms $\varphi,\psi$ of $D_{n,1}(\mu)$ and a positive integer $k$ such that
  $$\varphi\circ F\circ \psi(z_1,\cdots,z_n,w)=(\sqrt{k}z_1,\cdots,\sqrt{k}z_n,w^k).$$
\end{thm}

\vskip 6pt
In 1978, using the Cartan-Chern-Moser theory, Webster \cite{Web} took up again the problem of considering a proper holomorphic mapping $f$ from the $n$-ball $\mathbb{B}^n=\{z\in\mathbb{C}^n: |z|<1\}$ into the $(n+1)$-ball $\mathbb{B}^{n+1}=\{z\in\mathbb{C}^{n+1}: |z|<1\}$ and showed that there exist automorphisms $\sigma\in {\rm Aut}(\mathbb{B}^n)$ and $\tau\in {\rm Aut}(\mathbb{B}^{n+1})$ such that $\tau\circ f\circ \sigma=({\rm id},0)$, where $f$ is $C^3$-smooth up to the boundary and  $n>2$. In a subsequent paper, Faran \cite{Far} classified proper holomorphic mappings form $\mathbb{B}^2$ into $\mathbb{B}^3$ which are three times continuously differentiable up to the boundary. In another work, Cima-Suffridge \cite{Cima} studied certain reflection principle for CR mappings between hypersurfaces with codimension one and established the results of Webster and Faran for proper holomorphic mappings which are only twice continuously differentiable up to the boundary. In the same paper, Cima-Suffridge conjectured that any proper holomorphic mapping from $\mathbb{B}^n$ into $\mathbb{B}^N$ ($n>1$), which is $C^2$-smooth up to the boundary, must have the form $\tau\circ f\circ \sigma=(id,0)$ with automorphisms $\sigma\in {\rm Aut}(\mathbb{B}^n)$ and $\tau\in {\rm Aut}(\mathbb{B}^{N})$ when $N<2n-1$. Huang \cite{H} confirmed Cima-Suffridfe's conjecture and proved the following:

\vskip 6pt
\noindent
{\bf Theorem 1.E} (Huang \cite{H}) {\it  Let $M_1$ and $M_2$ be two connected open pieces of the boundaries of $\mathbb{B}^n\subset\mathbb{C}^n$ and $\mathbb{B}^N\subset \mathbb{C}^N$, respectively. Let $f$ be a nonconstant twice continuously differentiable CR mapping from $M_1$ into $M_2$. Suppose that $n>1, N<2n-1$. Then $f$ is the restriction of a certain totally geodesic embedding from $\mathbb{B}^n$ into $\mathbb{B}^N$. More precisely, there exist an automorphism $\sigma\in Aut(\mathbb{B}^n)$ and an automorphism $\tau\in Aut(\mathbb{B}^N)$ such that $\tau\circ f\circ \sigma(z_1,\cdots,z_n)\equiv (z_1,\cdots,z_n,0,\cdots,0)$.  }

\vskip 6pt
As a corollary, we have the result as follows.

\vskip 6pt
\noindent
{\bf Theorem 1.F} (Huang \cite{H}) {\it Let $f$ be a proper holomorphic mapping from $\mathbb{B}^n$ into $\mathbb{B}^N$, which is twice continuously differentiable up to the boundary. Suppose that $n>1, N<2n-1$. Then there exist $\sigma\in Aut(\mathbb{B}^n)$ and $\tau\in Aut(\mathbb{B}^N)$ such that $\tau\circ f\circ \sigma(z_1,\cdots,z_n)\equiv (z_1,\cdots,z_n,0,\cdots,0)$.  }

\vskip 6pt
Using this result, we obtain a classification of proper holomorphic  mappings between two nonequidimensional Fock-Bargman-Hartogs domains as follows:

\begin{thm}
  Let $D_{n,1}(\mu)$ and $D_{N,1}(\mu)$ be two Fock-Bargmann-Hartogs domains of dimension $n+1$ and $N+1$, respectively. Let $F$ be a proper holomorphic mapping from $D_{n,1}(\mu)$ into $D_{N,1}(\mu)$ that is twice continuously differentiable up to the boundary. Suppose that $N<2n$.
  Then there exist automorphisms $\sigma\in Aut(D_{n,1}(\mu))$, $\tau\in Aut(D_{N,1}(\mu))$ and a positive integer $k$ such that $\tau\circ F\circ \sigma(z_1,\cdots,z_n,w)\equiv (\sqrt{k}z_1,\cdots,\sqrt{k}z_n,0,\cdots,0,w^k)$.
\end{thm}

\section{The automorphism group of the unit ball}
Let $H_{p,q}$ denote the indefinite Hermitian bilinear form on $\mathbb{C}^{p+q}, p,q\geq 1$, defined by $H_{p,q}(w,w)=\sum_{1\leq i\leq p} |w_i|^2-\sum_{p+1\leq j\leq p+q} |w_j|^2$. Denote the group of linear isometries of $H_{p,q}$ by $U(p,q)\subset GL(p+q,\mathbb{C})$ and call it the unitary group of (the indefinite form) $H_{p,q}$. Those of determinant $1$ constitute the special unitary group $SU(p,q)$ of $H_{p,q}$.

Denote the unit ball by $\mathbb{B}^n:=\{ (z_1,\cdots,z_n)\in\mathbb{C}^n : \sum_{i=1}^n |z_i|^2 <1\}$. By embedding $\mathbb{C}^n$ as part of $\mathbb{P}^n$ via $(z_1,\cdots,z_n)\rightarrow [1,z_1,\cdots,z_n]$ in terms of homogeneous coordinates, we have $\mathbb{B}^n=\{[w_0,\cdots,w_n]\in\mathbb{P}^n: |w_0|^2-\sum_{1\leq i\leq n}|w_i|^2 >0\}$. The homogeneity of $\mathbb{B}^n$ can be seen by considering the automorphisms
$$\Psi(z_1,\cdots,z_n)=\Big(\frac{z_1+\alpha}{1+\overline{\alpha}z_1},\frac{\sqrt{1-|\alpha|^2}}{1+\overline{\alpha}z_1}z_2,\cdots, \frac{\sqrt{1-|\alpha|^2}}{1+\overline{\alpha}z_1}z_n \Big),$$
where $|\alpha|<1$. Given any $x\in\mathbb{B}^n$ there exists an automorphism $\Phi_x$ such that $\Phi_x(0)=x$ of the form $\Phi_x(z)=\Psi(U(z))$, where $U$ is a unitary transformation and $\Psi$ is of the form above for some choice of $\alpha$.

It is clear that $SU(1,n)$, acting as projective linear transformations on $\mathbb{P}^n$, preserves $\mathbb{B}^n$.
Conversely, the above transformations $\Psi$ are represented in homogeneous coordinates by
$$\Psi([w_0,\cdots,w_n])=\big[ w_0+\bar{\alpha}w_1,\alpha w_0+w_1,\sqrt{1-|\alpha|^2} w_2,\cdots, \sqrt{1-|\alpha|^2}w_n \big].$$
Writing $\beta=\sqrt{1-|\alpha|^2}$. Then $\Psi([w_0,\cdots,w_n])= \big[\frac{w_0+\bar{\alpha}w_1}{\beta},\frac{\alpha w_0+w_1}{\beta}, w_2,\cdots,w_n \big].$
Thus $\Psi$ can be represented by the following matrix
\[ \left(\begin{array}{ccc}
1/\beta & \bar{\alpha}/\beta & 0 \\
\alpha/\beta & 1/\beta  & 0 \\
0       &   0       & I_{n-1}
\end{array}\right)\in SU(1,n).
\]
We have more precisely $Aut(\mathbb{B}^n)\cong SU(1,n)/\{\mu_{n+1} I\}$ as a consequence of Cartan's theorem,
where $\{\mu_{n+1} I\}$ consisting of diagonal matrices $\varepsilon I$, where $\varepsilon$ is an $(n+1)$-th root of unity and $I=I_{n+1}$ is the identity matrix of rank $n+1$.

\begin{lem}
  Let $\Psi\in Aut(\mathbb{B}^{n+1})$ be a holomorphic automorphism of the unit ball $\mathbb{B}^{n+1}=\big\{(z_1,\cdots,\\ z_n, z_{n+1})\in\mathbb{C}^{n+1}\big|\sum_{1\leq i\leq n+1}|z_i|^2 < 1\big\}$ such that $\Psi(0,\cdots,0,1)=(0,\cdots,0,1)$. Then, up to a unitary transformation in $z_1,z_2,\cdots,z_n$, $\Psi=(\psi_1,\cdots,\psi_n,\psi_{n+1})$ has the following form:
  \begin{align*}
    &\psi_1(z_1,\cdots,z_{n+1})=\frac{a_1+z_1-a_1z_{n+1}}{a_0+\lambda\overline{a_1}z_1+\cdots+\lambda\overline{a_n}z_n+(\lambda-a_0)z_{n+1}},\\
    &\quad\quad\quad \vdots\\
    &\psi_n(z_1,\cdots,z_{n+1})=\frac{a_1+z_n-a_nz_{n+1}}{a_0+\lambda\overline{a_1}z_1+\cdots+\lambda\overline{a_n}z_n+(\lambda-a_0)z_{n+1}},\\
    &\psi_{n+1}(z_1,\cdots,z_{n+1})=\frac{\frac{1}{\lambda}+a_0+\lambda\overline{a_1}z_1+\cdots+\lambda\overline{a_n}z_n+(\lambda-\frac{1}{\lambda}-a_0)z_{n+1}}{a_0+\lambda\overline{a_1}z_1+
    \cdots+\lambda\overline{a_n}z_n+(\lambda-a_0)z_{n+1}},
  \end{align*}
  where $\lambda$ is a pure imaginary number.
  In the other words, the special unitary matrix corresponding to $\Psi$ has the following form:
  \[\left(\begin{array}{ccccc}
  a_0 & \lambda\overline{a_1} & \cdots & \lambda\overline{a_n} & \lambda-a_0 \\
  a_1 & 1              & \cdots & 0              & -a_1\\
  \vdots &             & \ddots &                & \vdots\\
  a_n & 0 &\cdots & 1 & -a_n\\
  \frac{1}{\lambda}+a_0 & \lambda\overline{a_1} & \cdots & \lambda\overline{a_n} & \lambda-\frac{1}{\lambda}-a_0
  \end{array}\right).
  \]
\end{lem}
\noindent
{\bf Remark}. By the discussion before Lemma 2.1, $\Psi\in Aut(\mathbb{B}^{n+1})\cong SU(1,n)/\{\mu_{n+1} I\}$ is a linear fractional transformation which is smooth on the boundary $\partial\mathbb{B}^{n+1}$ of $\mathbb{B}^{n+1}$. Thus the condition ``$\Psi(0,\cdots,0,1)=(0,\cdots,0,1)$" in the Lemma 2.1 is well defined.

\noindent
{\bf Proof}. Let $U\in SU(1,n+1)$ be the special unitary matrix corresponding to $\Psi$ as follows:
\[U=\left(\begin{array}{cccc}
a_{00} & a_{01} & \cdots & a_{0\;n+1}\\
a_{10} & a_{11} & \cdots & a_{1\;n+1}\\
\vdots & \vdots & \vdots & \vdots\\
a_{n+1\;0}&a_{n+1\;1} &\cdots &a_{n+1\;n+1}
\end{array}\right).
\]
Since $U\in SU(1,n+1)$, we have
\[U\left(\begin{array}{cccc}
1 & & & \\
 & -1& &\\
 &   &\ddots& \\
 &   &     & -1
\end{array}\right) \overline{U}^t=\left(\begin{array}{cccc}
1 & & & \\
 & -1& &\\
 &   &\ddots& \\
 &   &     & -1
\end{array}\right).
\]
Thus
\begin{align}
  &|a_{00}|^2-\sum_{i=1}^{n+1} |a_{0i}|^2=1 ,\\
  &|a_{j0}|^2-\sum_{i=1}^{n+1} |a_{ji}|^2=-1, \quad for\; 1\leq j\leq n,\\
  & a_{j0}\overline{a_{k0}}-\sum_{i=1}^{n+1} a_{ji} \overline{a_{ki}}=0, \quad for\; j\neq k.
\end{align}
Since $\Psi(0,\cdots,0,1)=(0,\cdots,0,1)$, we have
\begin{align*}
  &\frac{a_{j0}+a_{j\;n+1}\cdot 1}{a_{00}+a_{0\;n+1}\cdot 1}=0 \quad (1\leq j\leq n),\\
  &\frac{a_{n+1\;0}+a_{n+1\;n+1}\cdot 1}{a_{00}+a_{0\;n+1}\cdot 1}=1.
\end{align*}
Thus,
\begin{align}
  &a_{j\;n+1}=-a_{j\;0} \quad ( 1\leq j\leq n),\\
  &a_{n+1\;0}+a_{n+1\;n+1}=a_{00}+a_{0\;n+1}=:\lambda.
\end{align}
Hence,
\begin{align*}
  &|a_{j\;0}|^2=|a_{j\;n+1}|^2,\quad 1\leq j\leq n, \\
  & a_{j0}\overline{a_{k0}}-a_{j\;n+1}\overline{a_{k\;n+1}}=0.
\end{align*}
Therefore,
\begin{align*}
  &\sum_{i=1}^n |a_{ji}|^2 =1, \quad 1\leq j\leq n, \\
  &\sum_{i=1}^n a_{ji} \overline{a_{ki}}=0, \quad j\neq k.
\end{align*}
That means that if we set $A=(a_{ij})_{1\leq i,j\leq n}$, then $A$ is a unitary matrix of rank $n$.

Define $\mathcal{A}:\mathbb{C}^{n+1}\rightarrow \mathbb{C}^{n+1},
\left(\begin{array}{c}
z_1\\ \vdots \\ z_{n+1}\end{array}\right)
\mapsto \left(\begin{array}{cc}
A^{-1} & 0\\
0 & 1
\end{array}\right)\left(\begin{array}{c}
z_1\\ \vdots \\ z_{n+1}\end{array}\right)$. Then $\mathcal{A}$ is a unitary transformation of $z_1,\cdots,z_n$ and its corresponding special unitary matrix in $SU(1,n+1)$ is
$\left(\begin{array}{ccc}
1 & & \\
  & A^{-1} & \\
  &        & 1
\end{array}\right)$.
Thus $\Phi:=\mathcal{A}\circ \Psi\in Aut(\mathbb{B}^{n+1})$ satisfying $\Phi(0,\cdots,0,1)=(0,\cdots,0,1)$ and  its corresponding special unitary matrix in $SU(1,n+1)$ is
\begin{align*}
 V=& \left(\begin{array}{ccc}
1 &  &  \\
  & A^{-1} & \\
  &        & 1\\
\end{array}\right) U   \\
=& \left(\begin{array}{ccccc}
a_{00} & a_{01} & \cdots & a_{0\;n}& a_{0\;n+1}\\
a_{10} & 1     &         &         & a_{1\;n+1}\\
\vdots &       & \ddots  &         & \vdots    \\
a_{n0} &       &         &  1      & a_{n\;n+1}\\
a_{n+1\;0} & a_{n+1\;1} & \cdots & a_{n+1\;n} & a_{n+1\;n+1}
\end{array}\right)  \\
=:& \left(\begin{array}{ccccc}
a_{0} & b_1 & \cdots & b_n & d_0\\
a_{1} & 1     &         &         & d_1\\
\vdots &       & \ddots  &         & \vdots    \\
a_{n} &       &         &  1      & d_n\\
a_{n+1} & c_1 & \cdots & c_n & d_{n+1}
\end{array}\right).
\end{align*}
By (2.4) and (2.5), we have $d_0=\lambda-a_0, d_{n+1}=\lambda-a_{n+1}, d_i=-a_i, 1\leq i\leq n$. Therefore
\[V=\left(\begin{array}{ccccc}
a_{0} & b_1 & \cdots & b_n & \lambda-a_0\\
a_{1} & 1     &         &         & -a_1\\
\vdots &       & \ddots  &         & \vdots    \\
a_{n} &       &         &  1      & -a_n\\
a_{n+1} & c_1 & \cdots & c_n & \lambda-a_{n+1}
\end{array}\right).
\]
Since $V\left(\begin{array}{cccc}
1 & & & \\
  &-1 & &\\
  &   &\ddots & \\
  &   &       & -1
\end{array}\right) \overline{V}^t=\left(\begin{array}{cccc}
1 & & & \\
  &-1 & &\\
  &   &\ddots & \\
  &   &       & -1
\end{array}\right)$, we get
\begin{equation}
 \begin{cases}
  &a_0 \cdot \overline{a_i}-b_i\cdot 1+(-\lambda+a_0)\cdot (-\overline{a_i})=0,\\
  &a_{n+1} \cdot \overline{a_i}-c_i\cdot 1+ (-\lambda+a_{n+1})\cdot (-\overline{a_i})=0,\\
  &|a_0|^2-\sum^n_{i=1} |b_i|^2-|\lambda-a_0|^2=1, \\
  &|a_{n+1}|^2-\sum^n_{i=1}|c_i|^2-|\lambda-a_{n+1}|^2=-1.
 \end{cases}
\end{equation}
Thus
\begin{align}
  b_i=\lambda \overline{a_i},\;\;\; c_i=\lambda\overline{a_i}.
\end{align}
Therefore
\[V=\left(\begin{array}{ccccc}
a_{0} & \lambda \overline{a_1} & \cdots & \lambda \overline{a_n} & \lambda-a_0\\
a_{1} & 1     &         &         & -a_1\\
\vdots &       & \ddots  &         & \vdots    \\
a_{n} &       &         &  1      & -a_n\\
a_{n+1} & \lambda \overline{a_1} & \cdots & \lambda \overline{a_n} & \lambda-a_{n+1}
\end{array}\right).
\]
Since $\det V=\lambda(a_{n+1}-a_0)=1$, we get $a_{n+1}-a_0=\frac{1}{\lambda}$. Thus
\[V=\left(\begin{array}{ccccc}
a_{0} & \lambda \overline{a_1} & \cdots & \lambda \overline{a_n} & \lambda-a_0\\
a_{1} & 1     &         &         & -a_1\\
\vdots &       & \ddots  &         & \vdots    \\
a_{n} &       &         &  1      & -a_n\\
\frac{1}{\lambda}+a_0 & \lambda \overline{a_1} & \cdots & \lambda \overline{a_n} & \lambda-\frac{1}{\lambda}-a_0
\end{array}\right).
\]
By (2.6) and (2.7), we get
$\lambda(\overline{a_0}-\overline{a_{n+1}})+\overline{\lambda}(a_0-a_{n+1})=2$. Therefore we have $-\frac{\lambda}{\overline{\lambda}}-\frac{\overline{\lambda}}{\lambda}=2$,
that is, $(\lambda+\overline{\lambda})^2=0$.  Hence $\lambda$ is a pure imaginary number. The proof of Lemma 2.1 is completed.

\begin{lem}
 Let $M_1$ and $M_2$ be two connected open pieces of the boundaries of $\mathbb{B}^n\subset\mathbb{C}^n$ and $\mathbb{B}^N\subset \mathbb{C}^N$, respectively, which contain the points $P=(0,\cdots,0,1)\in\partial \mathbb{B}^n$ and $Q=(0,\cdots,0,1)\in\partial \mathbb{B}^N$, respectively. Let $f$ be a nonconstant twice continuously differentiable CR mapping from $M_1$ into $M_2$ such that $f(P)=Q$. Suppose that $n>1, N<2n-1$. Then there exist an automorphism $\sigma\in Aut(\mathbb{B}^n)$ and an automorphism $\tau\in Aut(\mathbb{B}^N)$ which have the form in Lemma 2.1 such that $\tau\circ f\circ \sigma(z_1,\cdots,z_n)\equiv (z_1,\cdots,z_n,0,\cdots,0)$.
\end{lem}

\noindent
{\bf Proof}. By Theorem 1.E, there exist $\phi\in Aut(\mathbb{B}^n)$ and $\psi\in Aut(\mathbb{B}^N)$ such that
$$\psi\circ f\circ \phi(z_1,\cdots,z_n)=(z_1,\cdots,z_n,0,\cdots,0)=(Id,0_{N-n})(z_1,\cdots, z_n).$$
So $f=\psi^{-1} \circ (Id,0_{N-n})\circ \phi^{-1}$. Let $P'=\phi^{-1}(P)$. Since $Aut(\mathbb{B}^n)$ is transitive on $\partial\mathbb{B}^n$, there exists $\rho\in Aut(\mathbb{B}^n)$ such that $\rho(P')=P$. Thus $\rho\circ \phi^{-1}(P)=P$. Therefore, if we let $\sigma=\phi\circ \rho^{-1}$, then $\sigma\in Aut(\mathbb{B}^n)$ and $\sigma(P)=P$. Hence, by Lemma 2.1, $\sigma$ has the required form.

Let $\varphi(z_1,\cdots,z_n,z_{n+1},\cdots,z_N)=(\rho^{-1}(z_1,\cdots,z_n), z_{n+1},\cdots,z_N)$. Then $\varphi\in Aut(\mathbb{B}^N)$ and
\begin{align*}
f(z_1,\cdots,z_n)=&\psi^{-1}\circ(Id,0_{N-n})\circ \phi^{-1}(z_1,\cdots,z_n)\\
&=\psi^{-1}(\phi^{-1}(z_1,\cdots,z_n),0_{N-n})\\
&=\psi^{-1}\Big(\rho^{-1}\big(\rho(\phi^{-1}(z_1,\cdots,z_n))\big),0_{N-n}\Big)\\
&=\psi^{-1}\circ \varphi\Big(\rho(\phi^{-1}(z_1,\cdots,z_n)),0_{N-n}\Big)\\
&=\psi^{-1}\circ \varphi\circ (Id,0_{N-n})\circ \rho\circ \phi^{-1} (z_1,\cdots,z_n).
\end{align*}
That is, $f=\psi^{-1}\circ \varphi\circ (Id,0_{N-n})\circ \rho\circ \phi^{-1}$.  If we set $\tau=\varphi^{-1}\circ \psi$, then $f=\tau^{-1}\circ (Id,0_{N-n})\circ \sigma^{-1}$.
That is, $\tau \circ f \circ \sigma= (Id,0_{N-n})$.
By $f(P)=Q$ and $(id,0_{N-n})\circ \sigma^{-1} (P)=(Id,0_{N-n})(P)=Q$, we have $\tau^{-1}(Q)=Q$. Hence, $\tau(Q)=Q$.
Therefore, by Lemma 2.1, $\tau$ has the required form. The proof of Lemma 2.2 is finished.

\vskip 10pt
In order to prove our main results, we also need the following lemma about the boundary regularity of proper holomorphic mappings between Fock-Bargmann-Hartogs domains:

\begin{lem}[Theorem 2.5 in Tu-Wang \cite{TW}]
If $D_{n,m}(\mu)$ and $D_{n',m'}(\mu')$ are two equidimensional Fock-Bargmann-Hartogs domains and $f$ is a proper holomorphic mapping from   $D_{n,m}(\mu)$ into $D_{n',m'}(\mu')$, then $f$ extends to be holomorphic in a neighborhood of $\overline{D_{n,m}(\mu)}$.
\end{lem}

\section{Proof of main results }

\noindent
{\bf Proof of Theorem 1.1}. Let $F$ be a proper holomorphic self-mapping of $D_{n,1}(\mu)$.
By Lemma 2.3, $F$ extends holomorphically to a neighborhood of $\overline{D_{n,1}(\mu)}$. Since the zero locus of the complex Jacobian $J_F$ of $F$ is a proper analytic subset, there is an open neighborhood $U$ of some boundary point $P\in \partial D_{n,1}(\mu)$ such that $F$ is a biholomorphic map from $U\cap D_{n,1}(\mu)$ to $D_{n,1}(\mu)$ with $f(U\cap \partial D_{n,1}(\mu))\subset \partial D_{n,1}(\mu)$.
Since $Aut(D_{n,1}(\mu))$ is transitive on $\partial D_{n,1}(\mu)$ by Theorem 1.C,  we can assume that $P=(0,\cdots,0,1)\in \partial D_{n,1}(\mu)$ and $F(P)=P$.
We can take $U=(\Delta_{\varepsilon}\times\cdots\times\Delta_{\varepsilon})\times \Delta_{\varepsilon}(1)\subset \mathbb{C}^n\times \mathbb{C}$, where
$\Delta_{\varepsilon}=\big\{z\in \mathbb{C}\big||z|<\varepsilon\big\}$ and $\Delta_{\varepsilon}(1)=\big\{w\in\mathbb{C}\big||w-1|<\varepsilon\big\}$ for some sufficient small positive number $\varepsilon$.

Define
\begin{align}
\varphi:& U\rightarrow \varphi(U)=:V, \\
&(z_1,\cdots,z_n,w)\mapsto (z_1,\cdots,z_n,-2i\log w)=:(z_1,\cdots,z_n,W),
\end{align}
where $\log w$ denotes the principal branch of logarithm on $\Delta_{\varepsilon}(1)$. Then, $\varphi$ is a biholomorphism between the neighborhood $U$ of the boundary point $P$ of $D_{n,1}(\mu)$ and the neighborhood $V:=\varphi(U)$ of the boundary point $O=(0,\cdots,0,0)$ of the Siegel upper half-space $H=\{Im W> |z_1|^2+\cdots+|z_n|^2\}$, and the inverse $\varphi^{-1}$ of  $\varphi$ is
\begin{align}
  \varphi^{-1}: &\varphi(U)\rightarrow U\\
  & (z_1,\cdots,z_n,W)\mapsto (z_1,\cdots,z_n,e^{\frac{i}{2}W}).
\end{align}
Moreover, $\varphi(P)=O$ and $\varphi(U\cap \partial D_{n,1}(\mu))\subset V\cap \partial H$.

Define
\begin{align}
  \psi:& V\rightarrow \psi(V)=:\Omega,\\
  &(z_1,\cdots,z_n,W)\mapsto (\frac{2z_1}{W+i},\cdots,\frac{2z_n}{W+i},-\frac{W-i}{W+i})=:(\xi_1,\cdots,\xi_n,\eta).
\end{align}
Then $\psi$ is a biholomorphism between the neighborhood $V$ of the boundary point $O$ of $H$ and the neighborhood $\Omega:=\psi(V)$ of the boundary point $Q=(0,\cdots,0,1)$ of the unit ball $\mathbb{B}^{n+1}$ in $\mathbb{C}^{n+1}$, and the inverse $\psi^{-1}$ of $\psi$ is
\begin{align}
  \psi^{-1}: &\Omega \rightarrow V\\
  & (\xi_1,\cdots,\xi_n,\eta)\mapsto (i\frac{\xi_1}{1+\eta},\cdots, i\frac{\xi_n}{1+\eta},i\frac{1-\eta}{1+\eta}).
\end{align}
Moreover, $\psi(O)=Q$ and $\psi(V\cap \partial H)\subset \Omega\cap \partial \mathbb{B}^{n+1}$.
Therefore,
$\psi\circ \varphi(U\cap \partial D_{n,1}(\mu))\subset \Omega\cap \partial \mathbb{B}^{n+1}$.

Define $G:=\psi\circ\varphi\circ F\circ (\psi\circ\varphi)^{-1}:\Omega\rightarrow \Omega$.  Then $G$ is biholomorphic, and
$G(\Omega\cap \partial \mathbb{B}^{n+1}) \subset \Omega\cap \partial \mathbb{B}^{n+1}$ and $G(Q)=Q$.
By Theorem 1.A, $G$ extends to a biholomorphic automorphism of $\mathbb{B}^{n+1}$ and $G(Q)=Q$.
Therefore, by Lemma 2.1, we have
\begin{align*}
& G(\xi_1,\cdots,\xi_n,\eta)\\
=&\big(g_1(\xi_1,\cdots,\xi_n,\eta),\cdots, g_n(\xi_1,\cdots,\xi_n,\eta),g_{n+1}(\xi_1,\cdots,\xi_n,\eta)\big)\\
=&\Big(\frac{a_1+\xi_1-a_1 \eta}{a_0+\lambda\overline{a_1}\xi_1+\cdots+\lambda\overline{a_n}\xi_n+(\lambda-a_0)\eta}, \cdots, \frac{a_n+\xi_n-a_n \eta}{a_0+\lambda\overline{a_1}\xi_1+\cdots+\lambda\overline{a_n}\xi_n+(\lambda-a_0)\eta}, \\ &\frac{(\frac{1}{\lambda}+a_0)+\lambda\overline{a_1}\xi_1+\cdots+\lambda\overline{a_n}\xi_n+(\lambda-\frac{1}{\lambda}-a_0)\eta}
{a_0+\lambda\overline{a_1}\xi_1+\cdots+\lambda\overline{a_n}\xi_n+(\lambda-a_0)\eta}  \Big).
\end{align*}
Therefore,
\begin{align*}
  & F(z_1,\cdots,z_n,w)\\
  =&(\psi\circ\varphi)^{-1} \circ G\circ (\psi\circ\varphi)(z_1,\cdots,z_n,w)\\
  =&\Big( i\frac{g_1(\psi\circ\varphi(z_1,\cdots,z_n,w))}{1+g_{n+1}(\psi\circ\varphi(z_1,\cdots,z_n,w))},\cdots, i\frac{g_n(\psi\circ\varphi(z_1,\cdots,z_n,w))}{1+g_{n+1}(\psi\circ\varphi(z_1,\cdots,z_n,w))},\\
  &\exp\Big\{-\frac{1}{2}\frac{1-g_{n+1}(\psi\circ\varphi(z_1,\cdots,z_n,w))}{1+g_{n+1}(\psi\circ\varphi(z_1,\cdots,z_n,w))}\Big\} \Big)\\
  =&\Big( i\frac{ a_1+\frac{2z_1}{-2i\log w+i}+a_1\frac{-2i\log w-i}{-2i\log w+i} }{ \frac{1}{\lambda}+2a_0+ 2\lambda\sum_{k=1}^n \overline{a_k} \frac{2z_k}{-2 i\log w+i} +\big[2(\lambda-a_0)-\frac{1}{\lambda}\big]\big(-\frac{-2i\log w-i}{-2i\log w+i}\big)}
  \cdots, \\
  & i\frac{ a_n+\frac{2z_n}{-2i\log w+i}+a_n\frac{-2i\log w-i}{-2i\log w+i} }
  { \frac{1}{\lambda}+2a_0+ 2\lambda\sum_{k=1}^n \overline{a_k} \frac{2z_k}{-2 i\log w+i} +\big[2(\lambda-a_0)-\frac{1}{\lambda}\big]\big(-\frac{-2i\log w-i}{-2i\log w+i}\big)}, \\
  & \exp\Big\{
  -\frac{1}{2}\frac{-\frac{1}{\lambda}+  \frac{1}{\lambda} \big(-\frac{-2i\log w-i}{-2i\log w+i}\big)}
  { \frac{1}{\lambda}+2a_0+ 2\lambda\sum_{k=1}^n \overline{a_k} \frac{2z_k}{-2 i\log w+i} +\big[2(\lambda-a_0)-\frac{1}{\lambda}\big]\big(-\frac{-2i\log w-i}{-2i\log w+i}\big)}
  \Big\} \Big) \\
   =&\Big( i\frac{2z_1+ a_1(-4i\log w) } {4\lambda \sum_{k=1}^n \overline{a_k} z_k + 2i(2\lambda - 4a_0 -\frac{2}{\lambda})\log w + 2i\lambda},\cdots, \\
  & i\frac{2z_n + a_n(-4i\log w) } {4\lambda \sum_{k=1}^n \overline{a_k} z_k + 2i(2\lambda - 4a_0 -\frac{2}{\lambda})\log w + 2i\lambda}, \\
  & \exp\Big\{
  -\frac{1}{2}\frac{\frac{4}{\lambda}i\log w }
  {4\lambda \sum_{k=1}^n \overline{a_k} z_k + 2i(2\lambda - 4a_0 -\frac{2}{\lambda})\log w + 2i\lambda}
  \Big\} \Big).
\end{align*}
Since $F|_{\{z_1=\cdots=z_n=0\}\times \{|w|<1\}}$ is holomorphic, we have that
$\frac{a_k(-4i\log w) } {2i(2\lambda - 4a_0 -\frac{2}{\lambda})\log w + 2i\lambda}$ $(1\leq k\leq n)$ are holomorphic functions on the unit disc $\{|w|<1\}$. Hence, $a_k=0, 1\leq k\leq n$. Thus,
\begin{align*}
  F(z_1,\cdots,z_n,w)
  =&\Big( i\frac{2z_1 } { 2i(2\lambda - 4a_0 -\frac{2}{\lambda})\log w + 2i\lambda},\cdots,
   i\frac{2z_n } { 2i(2\lambda - 4a_0 -\frac{2}{\lambda})\log w + 2i\lambda}, \\
  & \exp\Big\{
  -\frac{1}{2}\frac{\frac{4}{\lambda}i\log w }
  {2i(2\lambda - 4a_0 -\frac{2}{\lambda})\log w + 2i\lambda}
  \Big\} \Big).
\end{align*}
Since $F|_{\{z_1=\varepsilon, z_2=\cdots=z_n=0\}\times \{|w|<\delta\}}$ is holomorphic, we have that
$\frac{2\varepsilon } { 2i(2\lambda - 4a_0 -\frac{2}{\lambda})\log w + 2i\lambda}$ is a holomorphic function on the small disc $\{|w|<\delta\}$.
Therefore,  $2\lambda - 4a_0 -\frac{2}{\lambda}=0$. Hence,
\begin{align*}
  F(z_1,\cdots,z_n,w)
  =&\Big( \frac{z_1 }{\lambda},\cdots, \frac{z_n}{\lambda},  \exp\Big\{-\frac{\log w }{\lambda^2}\Big\} \Big)
  =\Big( \frac{z_1 }{\lambda},\cdots, \frac{z_n}{\lambda},  w^{-\frac{1}{\lambda^2}} \Big).
\end{align*}
Since $\lambda$ is an pure imaginary number, we known that $-\frac{1}{\lambda^2}$ is a positive real number.
Moreover, since $F$ is holomorphic at $(0,\cdots,0,0)$, we have that $w^{-\frac{1}{\lambda^2}}$ is holomorphic at $w=0$. That is, $-\frac{1}{\lambda^2}\in \mathbb{Z}_+$ is a positive integer.
Therefore there exist holomorphic automorphisms $\varphi,\psi$ of $D_{n,1}(\mu)$ and a positive integer $k$ such that
  $$\varphi\circ F\circ \psi(z_1,\cdots,z_n,w)=(\sqrt{k}z_1,\cdots,\sqrt{k}z_n,w^k).$$
The  proof of Theorem 1.1 is completed.

\vskip 10pt
\noindent
{\bf Proof of Theorem 1.2}.
Since $Aut(D_{n,1}(\mu))$ is transitive on $\partial D_{n,1}(\mu)$ and $Aut(D_{N,1}(\mu))$ is transitive on $\partial D_{N,1}(\mu)$ by Theorem 1.C,   we can assume that $P=(0,\cdots,0,1)\in \partial D_{n,1}(\mu)$ and $Q=(0,\cdots,0,1)\in \partial D_{N,1}(\mu)$ and $F(P)=Q$. Since $F$ is twice continuously differentiable up to the boundary, there exists
$U=(\Delta_{\varepsilon}\times\cdots\times\Delta_{\varepsilon})\times \Delta_{\varepsilon}(1)\subset \mathbb{C}^n\times \mathbb{C}$ with $0<\varepsilon <1$ small enough such that $F$ is twice continuously differentiable on $U$.

Define
\begin{align*}
  \Phi: U &\rightarrow \Phi(U)=:V, \\
  (z_1,\cdots,z_n,w)&\mapsto \Big(\frac{2z_1}{-2i\log w+i}, \cdots, \frac{2z_n}{-2i\log w+i}, -\frac{-2i\log w-i}{-2i\log w+i} \Big)=:(\xi_1,\cdots,\xi_n,\eta),
\end{align*}
where $\log w$ denotes the principal branch of logarithm on $\Delta_{\varepsilon}(1)$. Then, $\Phi$ is a biholomorphism between the neighborhood $U$ of the boundary point $P$ of $D_{n,1}(\mu)$ and the neighborhood $V:=\Phi(U)$ of the boundary point $P'=(0,\cdots,0,1)$ of the unit ball $\mathbf{B}^{n+1}$, and the inverse $\Phi^{-1}$ of  $\Phi$ is
\begin{align*}
  \Phi^{-1}: V &\rightarrow U, \\
  (\xi_1,\cdots,\xi_n,\eta) &\mapsto \Big(i\frac{\xi_1}{1+\eta},\cdots, i\frac{\xi_n}{1+\eta}, \exp\big\{\frac{-1}{2}\frac{1-\eta}{1+\eta}\big\} \Big).
\end{align*}
Moreover, $\Phi(P)=P'$ and $\Phi(U\cap \partial D_{n,1}(\mu))\subset V\cap \partial \mathbf{B}^{n+1}$.

Similarly, let $W$ be a small open neighborhood of $Q$ and define
\begin{align*}
  \Psi: W &\rightarrow \Psi(W)=:\Omega, \\
  (Z_1,\cdots,Z_N,w)&\mapsto \Big(\frac{2Z_1}{-2i\log W+i}, \cdots, \frac{2Z_N}{-2i\log W+i}, -\frac{-2i\log W-i}{-2i\log W+i} \Big)=:(\zeta_1,\cdots,\zeta_N,\theta),
\end{align*}
where $\log w$ denotes the principal branch of logarithm on $\Delta_{\varepsilon}(1)$. Then, $\Psi$ is a biholomorphism between the neighborhood $W$ of the boundary point $Q$ of $D_{N,1}(\mu)$ and the neighborhood $\Omega:=\Psi(W)$ of the boundary point $Q'=(0,\cdots,0,1)$ of the unit ball $\mathbf{B}^{N+1}$, and the inverse $\Psi^{-1}$ of  $\Psi$ is
\begin{align*}
  \Psi^{-1}: V &\rightarrow U, \\
  (\zeta_1,\cdots,\zeta_N,\theta) &\mapsto \Big(i\frac{\zeta_1}{1+\theta},\cdots, i\frac{\zeta_N}{1+\theta}, \exp\big\{\frac{-1}{2}\frac{1-\theta}{1+\theta}\big\} \Big).
\end{align*}
Moreover, $\Psi(Q)=Q'$ and $\Psi(W\cap \partial D_{N,1}(\mu))\subset \Omega\cap \partial \mathbf{B}^{N+1}$.

Define $G:=\Psi\circ F\circ \Phi^{-1}:V\rightarrow \Omega$. Then $G$ is a nonconstant twice continuously differentiable CR mapping from $V\cap \partial\mathbf{B}^{n+1}$ into $\Omega\cap \partial\mathbf{B}^{N+1}$ such that $G(P')=Q'$. By Lemma 2.2, there exist $\sigma\in Aut(\mathbf{B}^{n+1})$ and $\tau\in Aut(\mathbf{B}^{N+1})$ which have the form in Lemma 2.1 such that $\tau\circ G\circ \sigma=(Id, 0_{N-n})$. Thus, $G=\tau^{-1}\circ (Id, 0_{N-n}) \circ \sigma$ and
\begin{align*}
  F&=\Psi^{-1}\circ G\circ \Phi \\
  &=\Psi^{-1}\circ \tau^{-1}\circ (Id, 0_{N-n}) \circ \sigma\circ \Phi  \\
  &=(\Psi^{-1}\circ \tau^{-1}\circ \Psi)\circ \Psi^{-1}\circ(Id,0_{N-n})\circ \Phi\circ (\Phi^{-1}\circ \sigma \circ \Phi).
\end{align*}
By the proof of Theorem 1.1, we have that there exist constants $a,b\in\mathbb{C}$ such that
\begin{align*}
(\Phi^{-1}\circ \sigma \circ \Phi)(z_1,\cdots,z_n,w)=(\sqrt{a}z_1,\cdots,\sqrt{a}z_n, w^a),\\
(\Psi^{-1}\circ \tau^{-1}\circ \Psi)(Z_1,\cdots,Z_N,W)=(\sqrt{b}Z_1,\cdots,\sqrt{b}Z_N, W^b).
\end{align*}
On the other hand,
\begin{align*}
  & \Psi^{-1}\circ(Id,0_{N-n})\circ \Phi(z_1,\cdots,z_n,w)\\
  =&\Psi^{-1} \Big(\frac{2z_1}{-2i\log w+i}, \cdots, \frac{2z_n}{-2i\log w+i}, 0_{N-n}, \frac{-2i\log w-i}{-2i\log w+i} \Big)\\
  =&(z_1,\cdots,z_n,0_{N-n},w).
\end{align*}
Thus,
\begin{align*}
 & F(z_1,\cdots,z_n,w) \\
 =&(\Psi^{-1}\circ \tau^{-1}\circ \Psi)\circ \Psi^{-1}\circ(Id,0_{N-n})\circ \Phi\circ (\Phi^{-1}\circ \sigma \circ \Phi)(z_1,\cdots,z_n,w)\\
 =&(\sqrt{ab}z_1,\cdots,\sqrt{ab}z_n,0_{N-n}, w^{ab}).
\end{align*}
Since $F$ is holomorphic at $(0,\cdots,0,0)$, we have that $ab \in\mathbb{Z}_+$. Thus, there exist automorphisms $\sigma\in Aut(D_{n,1}(\mu))$, $\tau\in Aut(D_{N,1}(\mu))$ and a positive integer $k$ such that
 $$\tau\circ F\circ \sigma(z_1,\cdots,z_n,w)\equiv (\sqrt{k}z_1,\cdots,\sqrt{k}z_n,0,\cdots,0,w^k).$$
The  proof of Theorem 1.2 is finished.

\vskip 10pt

\noindent\textbf{Acknowledgments}\quad  The authors would like to thank Professor Xianyu Zhou for his helpful disscussions, and thank the referees for useful comments. The first author was supported
by the National Natural Science Foundation of China (No. 11671306), and
the second author was partially supported by China Postdoctoral Science Foundation (No. 2016M601150).

\end{document}